\documentclass[a4paper]{IEEEtran}
\usepackage[latin1]{inputenc}
\usepackage[cmex10]{amsmath}
\usepackage{amsfonts}
\usepackage{amssymb}
\usepackage{graphicx}
\usepackage{verbatim}
\usepackage{array}
\usepackage{multirow}
\usepackage{dcolumn}
\usepackage{color}
\usepackage[noadjust]{cite}
\usepackage{url}
\usepackage{balance}
\usepackage[usenames,dvipsnames]{xcolor}

\usepackage{amsmath}

\usepackage{nomencl}
\usepackage{breqn}
\makenomenclature
\usepackage{etoolbox}

\DeclareGraphicsExtensions{.eps}

\begin{document}
\title{Addressing the Time-Varying Dynamic Probabilistic Reserve Sizing Method on Generation and Transmission Investment Planning Decisions}

\author{Alessandro~Soares,~Ricardo~Perez,~Weslly~Morais,~Silvio~Binato

\thanks{Alessandro Soares, Ricardo Perez, Weslly Morais and Silvio Binato are from PSR, Rio de Janeiro, RJ, Brazil (e-mail: \mbox{alessandro@psr-inc.com}; \mbox{ricardo@psr-inc.com}; \mbox{weslly@psr-inc.com}; \mbox{silvio@psr-inc.com})}
}

\maketitle


\begin{abstract}
In this paper, we address the long-term system's requirement reserve sizing due to the high-level of variable renewable energy (VRE) sources penetration, inside the expansion planning model. The increase in the insertion of this kind of energy source will also bring an increase in the reserve requirements. A higher requirement will be translated into additional costs to the system since the system operator will need to allocate generators for reserve purposes. The VRE sources implicitly cause these costs, so besides the investment cost, the expansion planning models should consider those costs on the expansion decision process. The methodology proposed here aims to provide a probabilistic and dynamic evaluation of the forecast errors of VRE sources generation, translating these errors into the system's requirement reserve. This evaluation is done inside the expansion planning optimization model, treating the reserve requirement as an endogenous variable. Finally, a real case study of the Mexican system is presented, so that we can analyze the results of the methodology and the impacts of considering the reserve requirement along with the expansion planning decision.
\end{abstract}


\label{sec:nomenclature}
\printnomenclature


\section{Introduction}
\label{sec:introduction}

\IEEEPARstart{T}{he} use of optimization models applied to the expansion of electrical systems is a practice adopted as a way to help planning agents to make decisions that bring more significant benefit to society, with the aim of meeting the load growth with the lowest possible investment and operating costs, maintaining, in contrast, the criteria of reliability, security of supply and also contemplating energy and environmental policies of governmental interests. In this way, investment planning optimization models are widely used to address large-scale Variable Renewable Energy (VRE) sources integration \cite{bird2012lessons,connolly2010review,ostergaard2009reviewing}.

The production of VRE sources is characterized by high volatility in a short period, as analyzed in the review in \cite{bevrani2010renewable}. Because of that, the high penetration of these sources in electrical systems has changed the way the systems are operated and planned. One consequence is that the tasks associated with the activities of operation and expansion planning of electrical systems must model these characteristics, inherent to variable renewable sources. One side effect of high VRE penetration is the increase of operational reserves requirements since dispatchable plants must have available capacity to compensate changes (increase or reductions) in the intermittent renewable generation to maintain the system reliability and stability. The work in \cite{eto2011use} analyzes in detail all the technical aspects for ensuring reliable operation in systems with a high amount of VRE penetration. Therefore, operational reserves should also have a dynamic behavior, varying within the hours of the day and seasons of the year, as they are associated with forecasting errors of VREs production, and should be based on a probabilistic process since they must model the stochastic process of the RES generation.


Since system's security and reliability is starting to be a problem for power system planning, the VRE uncertainty, along with other sources of system imbalance, increases the need for operational flexibility, i.e., the capacity to fast response for upward or downward generation ramping and to provide reserve capacity \cite{bouffard2011value,ulbig2012operational,nosair2015flexibility}. The existing conventional sources of flexibility is not enough for dealing with the recent penetration of VRE sources \cite{palmintier2015impact}. Because of that, there are some implicit flexibility costs attached to VRE sources that must be taken into account in the expansion planning decision. 

There are lots of conventional ways to address the system's security, such as using reliability indices \cite{hemmati2013reliability,kannan2005application}, and with capacity reserve margin target (the difference between installed capacity and the peak load) \cite{tafreshi2012reliable}. These approaches differ from the proposed Time-Varying Dynamic Probabilistic Reserve (TDPR), since they do not capture the physical impacts of the reserve in the load supply, and does not guarantee the reserve capacity availability when it is needed. Because of that, it is important to co-optimize energy and reserve allocations inside the expansion planning model, calculating the reserve requirement as an endogenous variable. 

There are several methodologies proposed for estimating operating reserve requirement due to wind uncertainty in the literature. \cite{soder1993reserve} proposed for the first time, in 1993, a methodology for estimating reserve requirement, using a normally distributed function to address the system's generation uncertainties. \cite{gouveia2004operational} evaluate the reserve requirement from the perspective of the operator, by quantifying the risk of not satisfying the net load (demand - VRE generation) using a monte-carlo simulation. \cite{garver1966effective,allan2013reliability} established stochastic methods for evaluating reserve requirement and \cite{matos2010setting} proposed similar methodologies but addressing the wind power forecast error (WPFE). \cite{bruninx2014statistical,holttinen2012methodologies,jost2015new} proposed that wind power uncertainty must not be treated as generator outages. These works proposed different distributions to model WPFE to estimate short-term requirement reserves. \cite{holttinen2012methodologies} proposed a gaussian distribution, \cite{bruninx2014statistical} proposed a Levy $\alpha$-stable distribution and \cite{jost2015new} proposed a non-parametric method named, Kernel Distribution.

An endogenous reserve requirement calculation has also been analyzed in the literature. The authors in \cite{van2017impact} proposed an expansion planning model with short-term operative constraints and an endogenous reserve calculation, assuming that reserve requirement has linear dependency with VRE generation. They analyzed in detail the impacts of considering an endogenous calculation of a dynamic reserve requirement inside the expansion planning model. The work in \cite{papavasiliou2011reserve} proposes a two-stage operating model, including Unit-Commitment (UC) and reserve balance constraints, considering the total system's requirement reserve as given. \cite{bucksteeg2016impacts} proposes a non-convex method for estimating this requirement, showing that a dynamic probabilistic reserve may have a substantial impact on the system's total cost.

In this paper, we propose an endogenous calculation of the system's requirement reserve, the Time-varying Dynamic Probabilistic Reserve (TDPR), for long-term horizon, due to the renewable intermittency inside a G\&T expansion planning optimization model. The methodology is applied to VRE generation scenarios, taking into account a statistical evaluation of forecast errors for each hour and VRE source, encompassing existing and candidate plants. Besides that, since the forecast error of the total VRE generation is considered, our methodology is capable of modeling the fact that complementary VRE sources may reduce the requirement reserve (a portfolio effect). The TDPR is formulated as a convex combination between the average and the Conditional Value at Risk (CVaR) of the probability distribution of the forecast error variation. Because of all these features, our model improves the methodology in \cite{van2017impact}, since the endogenous calculation proposed by the author does not consider portfolio effect and forecast errors. So, the proposed methodology is state-of-the-art in estimating requirement reserve for systems with a high level of VRE penetration.

The aforementioned expansion model is applied, considering short-term constraints, such as UC, ramping, VRE uncertainties, endogenously reserve commitment, and continuous/binary investment decision variables in a co-optimization scheme. i.e., the expansion model optimizes both expansion and reserves costs in an integrated framework. All of these details may result in a large, computationally intractable problem for real-world systems, especially because expansion planning problems usually have long horizons. In order to reduce computational effort, a time-clustering approach based on representative days is used. The resulting model is a Mixed-integer Linear Programming (MILP) that aims at minimizing investment costs and expected value of operating costs, where the operation is solved with hourly time steps.

The main contribution of this paper is a proposed method to endogenously estimate the total system's TDPR taking into account the regional portfolio effect achieved by the spatial correlation of the intermittent VRE sources. And as a consequence, co-optimize generation and transmission investment decisions jointly with TDPR.

In order to achieve this goal, a case study based on the Mexican system expansion will be analyzed, where investment decisions, technology mix, and total costs are compared when running a generation expansion planning with and without the proposed reserve formulation. Additionally, portfolio effects are analyzed by comparing (i) a hierarchical approach, where RES plants are decided without TDPR followed by investment decisions to accomplish the resulting necessary reserve requirement; with (ii) a co-optimization approach, where VRE sources are chosen with TDPR along with dispatchable plants to meet this requirement.

This paper is organized as follows. In Section \ref{sec:methodology}, we show the proposed methodology. In Section \ref{sec:case}, we perform a case study utilizing data from the Mexican power system in order to analyze the results of the methodology applied to a real system. Finally, in Section \ref{sec:conclusion}, the final conclusions are presented.


\section{Methodology}\label{sec:methodology}

The proposed methodology adds the endogenous reserve requirement calculation constraints to the expansion planning model proposed in \cite{optgen2}. Section \ref{sec:overall} will provide the overall methodology and Section \ref{sec:optm} how the methodology will be added as constraints in the optimization model for the endogenous calculation of the reserve.

\subsection{Overall methodology}\label{sec:overall}

Due to the intermittency and unpredictability of this kind of energy source, the calculation of the reserve requirement must be: (i) probabilistic, that is, take into account the stochastic process of VRE generation in consecutive hours; and (ii) dynamic, that is, taking into account the fact that VRE production varies throughout the hours of the day and throughout the months of the year. In practical terms, this means that the operational reserve due to VRE is represented as a time profile that varies by month (due to the seasonal pattern of production of the VRE) and by year (due to the entry of new capacity).

The TDPR computation is done by each month, and may be split into 4 steps. So, for each month, the steps below are used:
\begin{enumerate}
    \item Determine the forecast of the hourly generation profile of the VRE sources, as a profile of 24 hours (forecast of the generation during a typical day of the month). Several methods may be used to determine this forecast, since it uses hourly generation scenarios as an input data, this forecast is always an input for the TDPR calculation. There are several VRE generation forecast models in the literature that may be applied in the proposed methodology, such as \cite{bruninx2014statistical,holttinen2012methodologies,jost2015new,foley2012current,juban2007probabilistic,hoeltgebaum2018generating}. In order to simplify the current methodology, and since wind forecasting is not in the scope of this paper, we are using the average of the VRE hourly scenarios generation as the forecast model. The Equation \eqref{eq:forecast} is used in order to calculate the hourly forecast. For example, assuming there are 50 scenarios, and that each consists of 30 days x 24 hours / day = 720 hours of VRE production, we will have 50 x 30 = 1500 samples of VRE production for the first hour, the same for the second hour, and so on. The hourly production profile is the average of these 1500 values for each hour, as Equation \eqref{eq:forecast} shows.
    
    \nomenclature{$\hat{g}_{r,h}$}{Forecast generation of the VRE source $r$ at hour of the day $h$ [MW]}
    \nomenclature{$p_s$}{Probability of the scenario $s$}
    \nomenclature{$D$}{Number of days of the current month}
    
    \begin{equation}\label{eq:forecast}
        \hat{g}_{r,h} = \sum\limits_{s,d}\dfrac{p_s}{D}g_{r,s,d,h}
    \end{equation}
    
    \item Determine the errors in the forecast generation (unpredictable generation) - In this step we will compute the error of the total VRE generation of the system. In order to do so, the Equation \eqref{eq:error} sums the individual errors of all the VRE sources in the system, and since these errors may be upward or downward, some specific VRE sources errors may compensate other VRE sources errors. This occurs when there are sources with some kind of complementarity. Different from the step 1, the step 2 is not a input for the expansion planning model, since the amount of VRE sources in the system is a decision variable. Because of that, these errors may not be estimated through a probability distribution, since they are a decision variable of an optimization model.
    
    \nomenclature{$\delta_{s,d,h}$}{Forecast error of the system at the scenario $s$, day $d$ and hour of the day $h$ [MW]}
    
    \begin{equation}\label{eq:error}
        \delta_{s,d,h} = \sum\limits_{r}\left( g_{r,s,d,h} - \hat{g}_{r,h}\right)
    \end{equation}
    
    As an example of Equation \eqref{eq:error}, suppose that the VRE production at hour 1, for a specific scenario, is 9200 MW, and that the forecast at hour 1 is 9000 MW. In this case, we will have an error of 9200 - 9000 = 200 MW. These 200 MW corresponding to the "stochastic" component (unpredictable) of the VRE production, and which therefore requires a generation reserve for its management. The calculation of the error is repeated for each of the 1500 scenarios of hour 1; hour 2; etc. The final result is a matrix with 1320 lines (scenarios) and 24 columns (hours of the day). Each element of this matrix contains a deviation in MW, positive or negative, with respect to the average time profile.
    
    \item Determine the error variation of VRE production between consecutive hours, as showed in the Equation \eqref{eq:diff} - For example, suppose that the error for hour 1, scenario 1 is 200 MW and that for the next hour (hour 2, scenario 1), it is -300 MW (negative value). This means that there was an unpredictable reduction of -300-200 = -500 MW of the VRE production between hours 1 and 2. In turn, this points to the need for a upward reserve (that is, possibility of increasing the generation to compensate) of 500 MW for hour 1, scenario 1. This process is repeated for the 1500 scenarios of hours 1 and 2. 
    
    \nomenclature{$\Delta_{s,d,h}$}{Error forecast variation of the system at the scenario $s$, day $d$ and hour of the day $h$ [MW]}
    
    \begin{equation}\label{eq:diff}
        \Delta_{s,d,h} = \delta_{s,d,h-1} - \delta_{s,d,h}
    \end{equation}
    
    The variation of the error between consecutive hours represents the fact that system operators allocate reserve capacity for the reserve requirement for each hour, so the allocation necessary for the following hour must be the difference of the forecast error of the following hour from the forecast error of the current hour, since this error has already been taken into account. This fact is illustrated in Figure \ref{fig:errorvariation}.
    
    \begin{figure}
        \centering
        \includegraphics[scale=0.6]{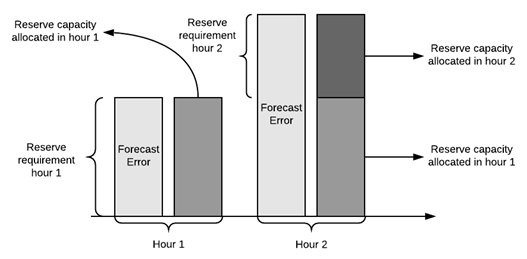}
        \caption{Reserve requirement calculation for hour 2}
        \label{fig:errorvariation}
    \end{figure}
    
    \item Determine the value of the upward reserve requirement of each hour, $TDPR_h$. Since most VRE sources can provide downward reserve \cite{milligan2010operating,joos2018short,van2017impact}, the proposed methodology only addresses upward reserve balance, since downward reserves are considered to be fully supplied by VRE sources. In order to do so, only the positive values of $\Delta_{s,d,h}$ should be considered, but this would cause the methodology to become non-convex, and would cause some problems while modelling inside the expansion planning optimization problem in Section \ref{sec:optm} (see \cite{burer2012non} for further details on non-convex optimization problems). Toward avoiding this issue, and considering that upward and downward reserve requirements have the same stochastic process, i.e, they have equal values when considering a big enough number of samples, rather than considering the positive or negative values of $\Delta_{s,d,h}$, the absolute value will be considered.
    
    The Equation \eqref{eq:tdpr} shows the TDPR calculation, that is a convex combination between the average and the Conditional Value of Risk (CVaR) \cite{rockafellar2000optimization}, so that it is possible to measure and decide the level of risks that system's planning is willing to take regarding system's reserve requirement.
    
    \nomenclature{$TDPR_h$}{The Time-Varying Dynamic Probabilistic Reserve at the hour of the day $h$ [MW]}
    \nomenclature{$\lambda$}{The risk-aversion factor}
    \nomenclature{$\beta$}{The reliability factor}
    \begin{equation}\label{eq:tdpr}
        TDPR_h = (1-\lambda)\sum\limits_{s,d}\dfrac{p_s}{D}|\Delta_{s,d,h}| + \lambda CVaR_{\beta}\left[ |\Delta_{s,d,h}| \right]
    \end{equation}
\end{enumerate}

\subsection{Inside the expansion planning optimization model}\label{sec:optm}
In order to calculate the TDPR inside the expansion planning model, the Equations \eqref{eq:forecast}, \eqref{eq:error}, \eqref{eq:diff} and \eqref{eq:tdpr} will be added as constraints to the optimization model described in \cite{optgen2}. The model \eqref{eq:modelbegin}-\eqref{eq:modelend} shows a simplified version of the model \cite{optgen2} and modifies the Equations \eqref{eq:forecast}, \eqref{eq:error}, \eqref{eq:diff} and \eqref{eq:tdpr} in order to add them as constraints:

\nomenclature{$g_{p,.}$}{Generation of the conventional energy source plant $p$ [MW]}
\nomenclature{$r_{p,.}$}{Reserve capacity of the conventional energy source plant $p$ [MW]}
\nomenclature{$g_{r,.}$}{Generation of the VRE source plant $r$ [MW]}
\nomenclature{$x_p$}{Investment decision variable of the conventional energy source plant $p$}
\nomenclature{$x_r$}{Investment decision variable of the VRE source plant $r$}
\nomenclature{$S$}{Total number of scenarios being considered}
\nomenclature{$c_p$}{Generating cost of the conventional energy source plant $p$ [$\$/MW$]}
\nomenclature{$I_p^T$}{Investment cost of the conventional energy source plant $p$ [$\$$]}
\nomenclature{$I_r^T$}{Investment cost of the VRE source plant $r$ [$\$$]}
\nomenclature{$g_{p,s,d,h}$}{Generation of the conventional energy source plant $p$ at scenario $s$, day $d$ and hour of the day $h$ [MW]}
\nomenclature{$g_{r,s,d,h}$}{Generation of the VRE source plant $r$ at scenario $s$, day $d$ and hour of the day $h$ [MW]}
\nomenclature{$\aleph$}{The set of investment constraint related to the expansion planning decision}
\nomenclature{$\tilde{G}$}{The set of operating constraint related to operational decisions}    \nomenclature{$\Delta_{s,d,h}^{+}$}{Positive value of the error forecast variation of the system at the scenario $s$, day $d$ and hour of the day $h$ [MW]}
\nomenclature{$\Delta_{s,d,h}^{-}$}{Negative value of the error forecast variation of the system at the scenario $s$, day $d$ and hour of the day $h$ [MW]}
\nomenclature{$W$}{Auxiliary variable used in the computation of the CVaR}
\nomenclature{$E_{h}$}{Average of the TDPR for the hour of the day $h$ [MW]}
\nomenclature{$CVaR_{h}$}{CVaR of the TDPR for the hour of the day $h$ [MW]}
\nomenclature{$\bar{g_p}$}{Maximum generation capacity of the conventional energy source plant $p$ [MW]}
    
\begin{align}
   &\min\limits_{g_{p,.},g_{r,.},x_p,x_r,r_{p,.}}  \ \dfrac{1}{S}\sum\limits_{p,s,d,h}c_pg_{p,s,d,h} + I_{p}^Tx_p + I_{r}^Tx_r\label{eq:modelbegin}\\
               s.t     &\ \ \ \ \{x_p,x_r\} \in \aleph\ & \label{eq:invcstr}\\
                       &\ \ \ \ \{g_{p,s,d,h},g_{r,s,d,h}\} \in \tilde{G}\ & \label{eq:opecstr} \\
                       &\ \ \ \ \delta_{s,d,h} = \sum\limits_{r}\left( g_{r,s,d,h} - \hat{g}_{r,h}\right)x_r  \label{eq:cooptimization}\\
                       &\ \ \ \ \Delta_{s,d,h}^{+} = \delta_{s,d,h-1} - \delta_{s,d,h}\label{eq:deltap}\\
                       &\ \ \ \ \Delta_{s,d,h}^{-} = -\left( \delta_{s,d,h-1} - \delta_{s,d,h}\right)\label{eq:deltan}\\
                       &\ \ \ \ \Delta_{s,d,h} = \Delta_{s,d,h}^{+}+\Delta_{s,d,h}^{-}\label{eq:delta}\\
                       &\ \ \ \ E_{h} = \sum\limits_{s,d}\dfrac{p_s}{D}\Delta_{s,d,h}\label{eq:average}\\
                       &\ \ \ \ CVaR_{h} = W+\dfrac{1}{\beta}\sum\limits_{s,d}p_s\omega_{s,d,h}\label{eq:cvar}\\
                       &\ \ \ \ TDPR_h \geq \left( 1-\lambda\right)E_{h}+\lambda CVaR_{h} \label{eq:tdprcstr}\\
                       &\ \ \ \ \omega_{s,d,h} + W \geq \Delta_{s,d,h}\label{eq:cvaraux}\\
                       &\ \ \ \ \omega_{s,d,h} \geq 0\\
                       &\ \ \ \ g_{p,s,d,h} + r_{p,s,d,h} \leq \bar{g_p}x_p\label{eq:reservealloc}\\
                       &\ \ \ \ \sum\limits_{p}r_{p,s,d,h} \geq TDPR_h\label{eq:reservecstr}\\
                       &\ \ \ \ \Delta_{s,d,h} \geq 0\\
                       &\ \ \ \ \Delta_{s,d,h}^{+} \geq 0\\
                       &\ \ \ \ \Delta_{s,d,h}^{-} \geq 0\label{eq:modelend}
\end{align}

The Constraint \eqref{eq:invcstr} represents all the investment constraints related to the model, such as capacity margin target, mix capacity margin, budget constraints and so on. Constraint \eqref{eq:opecstr} represents all operational constraints such as load balance, minimum and maximum generation, unit commitment constraints and so on. All the other constraints are related to the reserve co-optimization, where the methodology described in Section \ref{sec:overall} will be applied as constraints to the model.

The Constraint \eqref{eq:cooptimization} has been modified to consider VRE investment decision into the total forecast error. In order to do so, each individual forecast error is multiplied by its investment decision variables (we are considering that if a specific VRE source is existing, then the investment decision variable is equal 1) couples the VRE investment decision with the system reserve requirement, since the total forecast error will depends directly on this decision.

Constraints \eqref{eq:deltap}-\eqref{eq:delta} models the absolute value of the error forecast variation, represented in Equations \eqref{eq:diff} and \eqref{eq:tdpr}. The constraints \eqref{eq:average}-\eqref{eq:tdprcstr} models the convex combination, represented in Equation \ref{eq:tdpr}. The Constraints \eqref{eq:cvar} and \eqref{eq:cvaraux} models the CVaR of the error forecast absolute value.

The Constraints \eqref{eq:reservealloc} and \eqref{eq:reservecstr} are the reserve balance constraints. Constraint \eqref{eq:reservealloc} models the amount of reserve capacity available for each generator at each time period, that is, the difference between the maximum generation and the current generation is the amount of reserve capacity. Constraint \eqref{eq:reservecstr} ensure that the total reserve capacity is higher than the TDPR, and because of that, an increase in the TDPR, will also increase reserve capacity necessities and decrease total available generation capacity, which may lead to a solution where more expensive thermal resources needs to be dispatched in order to guarantee the load supply. So the Constraint \eqref{eq:reservecstr} model the implicit cost that VRE sources brings to the system through a higher reserve requirement.

\section{Case Study}\label{sec:case}
\subsection{Mexican system}\label{sec:mexicansystem}
\subsubsection{Main characteristics}
The National Electrical System (NES) has approximately 75 GW of installed capacity, where the main source is the natural gas thermal plants, with 54\% participation in the mix. Until 2017, wind plants constituted 5\% of the installed capacity, while solar plants constituted less than 1\%. The NES is divided into ten load control regions and 53 transmission regions.
\subsubsection{Expansion candidates}
The representation of the generators unit is based on the information presented in the PRODESEN database, where a short-term expansion plan, an indicative expansion and a decommissioning program, are presented. Candidates for expansion are presented in table \ref{tab:study}.

\begin{table}[]
\caption{Candidate projects and investment cost}
\label{tab:study}
\begin{tabular}{c|c|c|c|c|c|c}
                                                                    & \textbf{CCGT} & \textbf{OCGT} & \textbf{Diesel} & \textbf{Solar} & \textbf{Wind} & \textbf{Geo} \\ \hline
\begin{tabular}[c]{@{}c@{}}Investment\\ cost\\ (\$/kW)\end{tabular} & 700           & 600           & 700             & 900            & 1200          & 3800         \\ \hline
\begin{tabular}[c]{@{}c@{}}Lifetime\\ (years)\end{tabular}          & 20            & 20            & 20              & 25             & 25            & 25           \\ \hline
\begin{tabular}[c]{@{}c@{}}O\&M\\ Costs\\ (\$/kW year)\end{tabular} & 25            & 15            & 12              & 16             & 25            & 20           \\ \hline
\begin{tabular}[c]{@{}c@{}}Efficiency\\ (\%)\end{tabular}           & 56            & 41            & 45              &                &               &              \\ \hline
\end{tabular}
\end{table}

\begin{figure}
        \centering
        \includegraphics[scale=0.5]{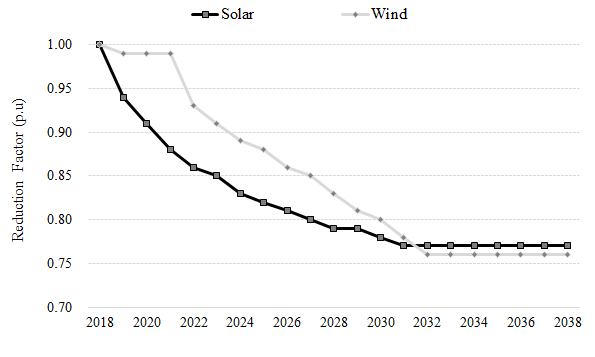}
        \caption{Reduction curve of the investment costs at Solar and Wind power plants}
        \label{fig:downcost}
    \end{figure}

As observed in recent years around the world, the investment cost of VRE sources follows a downward trend \cite{augustine20182018}. Thus the decay curve presented in Figure \ref{fig:downcost} was considered as a way to capture this trend. In order to simulate the stochastic hourly operation of these plants, a set of 32 hourly scenarios was considered, generated from the MERRA-2 database \cite{gelaro2017modern}.

\subsection{Results}
The case study presented here has the goal to analyze and study the impacts of the TDPR constraints in the expansion plan of a system with a large insertion of VRE sources in the energy matrix. The Mexican system will be simulated for the year 2038 using the configuration of 2018 as a starting point.

\begin{figure}
        \centering
        \includegraphics[scale=0.5]{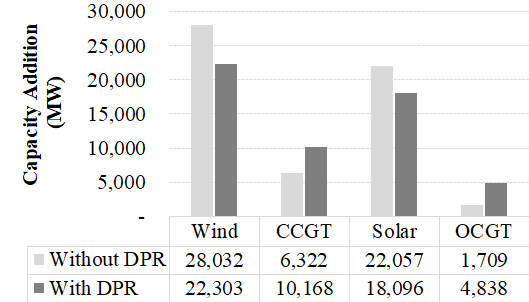}
        \caption{Capacity addition per technology type}
        \label{fig:results}
    \end{figure}

In this study, two expansion alternatives are analyzed: (i) expansion plan without TDPR constraints; and (ii) with TDPR constraints. Figure \ref{fig:results} shows the addition of installed capacity to the system by technology for the two expansion alternatives.

The case (i) "Without TDPR" can be defined as a "full economical alternative," since the model will minimize costs without considering the operational reserve constraints, i.e., disregarding the implicit cost of VRE sources. So the results of this alternative will be fully based on the trade-off between the investment cost and the reduction of the operating costs. 

The case (ii) "With TDPR" has a reduction in the added capacity of VRE sources and an increase of natural gas plants. Therefore, it is possible to conclude that this increase in dispatchable plants in the system is mainly because the TDPR is being considered inside the expansion planning process. There is an increase of 62\% for combined cycle plants (CCGT) and 183\% for open-cycle plants (OCGT). Although OCGTs are less efficient than combined-cycle plants (higher levelized cost of energy), the model decided to invest in this technology to meet reserve requirements, since these plants have lower investment costs. In general, OCGTs plants are built to meet peak demand and to operate in critical situations. Because of its higher operative costs, this kind of plant is not built to dispatch frequently, otherwise, the CCGTs would be better alternatives.

\begin{figure}
    \centering
    \includegraphics[scale=0.53]{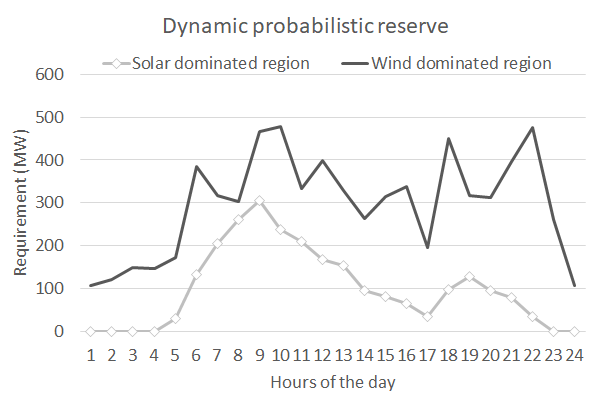}
    \caption{DPR for the two Mexican regions}
    \label{fig:rpdresults}
\end{figure}

Figure \ref{fig:rpdresults} shows the total TDPR requirement for two different regions of the Mexican system, one dominated by solar plants and the other one by wind plants. As one can see, the solar dominated area has two peaks in the requirement during the sunrise and the sunset. At night the requirement is zero since solar plants are not generating, and during the daylight, the requirement is small because the generation is usually flat with a high predictability level. 

Analyzing the wind dominated region, since the variability between hours of wind plants is usually higher than solar plants, the reserve requirement is also higher. For regions dominated by wind power plants, the profile is linked to the wind pattern of the region, thus it does not have a standard profile like the solar dominated regions.

\begin{figure}
    \centering
    \includegraphics[scale=0.6]{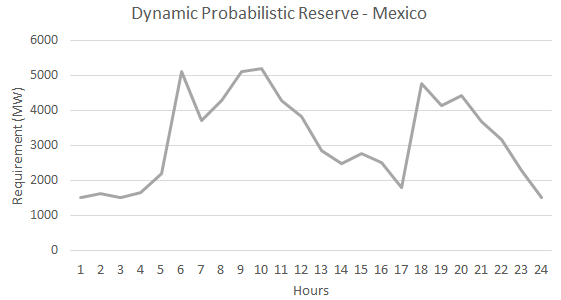}
    \caption{Total DPR for the Mexican system}
    \label{fig:totrpd}
\end{figure}

The Figure \ref{fig:totrpd} shows the total TDPR requirement for the Mexican system. As one can see, the maximum requirement is 5 GW, which is 12.5\% of the total installed capacity of VRE sources. Besides that, the TDPR tends to be bigger during the sunrise and the sunset, due to solar power plants, and to be smaller during the night.

\section{Conclusions}\label{sec:conclusion}
The case study showed that the co-optimization between generation/transmission planning and operational reserves could improve reliability in the supply of demand in a system with a massive insertion of VRE sources. The TDPR constraints inside the planning model made the expansion plan to be adapted to accommodate the high variability of VRE generation in the system by introducing natural gas power plants, increasing system's flexibility.

The results also showed that the model reduces VRE source capacity added to the system when compared to the case without the co-optimization with the TDPR. This reduction may happen because VRE sources variability may have an associated (implicit) cost to the system due to an increase in the reserve requirement. This increase is reflected in terms of total cost.

Additionally, it redistributes the VRE sources geographically to reduce the reserve requirement, due to the portfolio effect that may happen. It was also observed changes in the expansion of the transmission system, because, during the co-optimization process, dispatchable plants were geographically reallocated.

\bibliographystyle{IEEEtran}
\bibliography{Bibliography}

\begin{IEEEbiography}[{\includegraphics[width=1in,height=1.25in,clip,keepaspectratio]{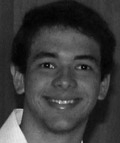}}]%
{Alessandro Soares}
Has a BSc in Electrical Engineering and in Control Engineering from PUC-Rio. Is currently doing a MSc in Optimization/Operations Research at PUC-RJ. Joined PSR in 2017 and has been working on the development and support of the expansion planning model (OPTGEN) and with the time series analysis model Time Series Lab (TSL). Before joining PSR, worked with non-linear time series models for synthetic inflow scenarios generation and modelling probabilistic scenarios for renewable plants. He has also worked with new methodologies for the analysis of the impact of climate changes in the inflow series.
\end{IEEEbiography}

\begin{IEEEbiography}[{\includegraphics[width=1in,height=1.25in,clip,keepaspectratio]{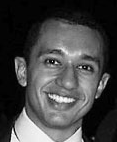}}]%
{Weslly Morais}
Electrical Engineer from UFRJ in 2018, with emphasis on Power Systems. Between 2014 and 2015, studied at the University of Southampton, UK, during an exchange program and worked at TE Connectivity, where he had the opportunity to be part of the team responsible for the development of products to be used under high voltage operation conditions (electric cables, circuit breakers and other protection systems). Joined PSR in 2016 and since then has been working in planning studies for the expansion of generation and transmission in several Latin American countries and in prices forecast studies and planning for the Brazilian electric and energy sector. His area of activity is the application of optimization models in electrical systems.
\end{IEEEbiography}

\begin{IEEEbiography}[{\includegraphics[width=1in,height=1.25in,clip,keepaspectratio]{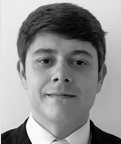}}]%
{Ricardo Perez}
Has a BSc in Electrical Engineering with emphasis in Power Systems from Itajuba Federal University (UNIFEI). During college he worked in two R\&D; projects in the Power Quality Study Group. Through an exchange program with the Technische Universitat Dresden in Germany, he also developed a research in the same field at this university. In addition to the applied research, his experience in the electricity business includes internships in Brazil and in Germany, respectively at the Generation Planning Department of CPFL Geracao and at DIgSILENT GmbH (where he carried out studies regarding the connection of Wind Farms to the German grid). Mr. Perez joined PSR in December 2009 and has been a member of the transmission studies group ever since.
\end{IEEEbiography}

\begin{IEEEbiography}[{\includegraphics[width=1in,height=1.25in,clip,keepaspectratio]{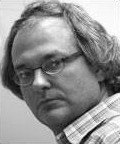}}]%
{Silvio Binato}
Is responsible for the development of methodologies and software for SDDP, OptGen and ePSR. He has a BSc and a MSc degrees in EE and a DSc in Systems Engineering and Computer Science. Dr. Binato has participated in expansion planning studies in Latin America and Europe. Previously, Dr. Binato worked at Cepel, where he managed R\&D; projects in the areas of transmission network expansion planning, optimal power flow and meta heuristics for global optimization in power systems. Dr. Binato is author or co-author of 25 papers in refereed journals and international conference proceedings.
\end{IEEEbiography}

\end{document}